\documentclass[12pt]{article}

\usepackage{amssymb}
\usepackage{amsmath}

\newtheorem{thrm}{Theorem}
\newtheorem{cor}{Corollary}
\newtheorem{lem}{Lemma}

\begin{document}

\title{Finitary incidence algebras}

\author{N.\,S.\,Khripchenko, B.\,V.\,Novikov}
\date{\footnotesize{Dept of Mechanics and Mathematics, Kharkov National University,
Ukraine\\
e-mail: boris.v.novikov@univer.kharkov.ua} }
\maketitle

It is well known, that incidence algebras can be defined only for locally finite
partially ordered sets~\cite{d-s-r, sta}. At the same time, for example, the poset of
cells of a noncompact cell partition of a topological space is not locally finite. On the
other hand, some operations, such as the order sum and the order product~\cite{sta}, do
not save the locally finiteness. So it is natural to try to generalize the concept of
incidence algebra.

In this article we consider the functions in two variables on an
arbitrary poset (\emph{finitary series}), for which the convolution
operation is defined. We obtain the generalization of incidence
algebra --- \emph{finitary incidence algebra} and describe its
properties: invertibility, the Jackobson radical, idempotents,
regular elements. As a consequence a positive solution of the
isomorphism problem for such algebras is obtained.

\section{Definition of finitary algebras}\label{sec-1}

In what follows $K$ denotes a fixed field, $P$ an arbitrary partially ordered set
(poset). Denote by $I(P)$ the set of formal sums of the form
\begin{equation}\label{eq-1}
\alpha = \sum_{x\le y}\alpha(x,y)[x,y]
\end{equation}
where $x,y\in P$, $\alpha(x,y)\in K$, $[x,y]= \{z\in P\mid x\le
z\le y\}$. In the general case $I(P)$ is not an algebra, so we
call it {\it an incidence space}.

A formal sum~\eqref{eq-1} is called {\it a finitary series} if for any $x,y\in P$, $x<y$,
there is only a finite number of subsegments $[u,v]\subset [x,y]$ such that $u\ne v$ and
$\alpha(u,v)\ne 0$. The set of all finitary series is denoted by $FI(P)$.

If $P$ is not locally finite then, as we mentioned above, the usual multiplication
(convolution) of the series
$$
\alpha\beta = \sum_{x\le y}\alpha(x,y)[x,y]\cdot\sum_{u\le
v}\beta(u,v)[u,v] = \sum_{x\le y}\left(\sum_{x\le z\le
y}\alpha(x,z)\beta(z,y)\right)[x,y]
$$
is not always defined in $I(P)$. The following proposition shows that in the general
situation it is reasonable to consider $FI(P)$ instead of $I(P)$.

\begin{thrm}\label{thrm-1} $FI(P)$ is an associative algebra and $I(P)$ is a module over it.
\end{thrm}

{\bf Proof.} Obviously, it is sufficient to prove that if $\alpha\in FI(P)$, $\beta\in
I(P)$ then $\alpha\beta$ is defined, and if in addition $\beta\in FI(P)$ then
$\alpha\beta\in FI(P)$.

So let $\alpha\in FI(P)$, $\beta\in I(P)$, $x,y\in P$ and $x\leq y$. Then the sum
$$
\gamma(x,y)=\sum_{x\le z\le y}\alpha(x,z)\beta(z,y)
$$
contains only a finite number of nonzero values $\alpha(x,z)$. Since every $\alpha(x,z)$
appears in the sum at most one time, this sum is finite. Therefore $\alpha\beta$ is
defined.

Now let $\alpha,\beta\in FI(P)$. Suppose that $\gamma=\alpha\beta\not\in FI(P)$. This
means that we can find such $x,y\in P$, $x\leq y$, that there is an infinite number of
subsegments $[u_i,v_i]\subset [x,y]$ $(i=1,2,\ldots)$, for which $u_i\ne v_i$ and
$\gamma(u_i,v_i)\ne 0$. At least one of the sets $\{u_i\}$, $\{v_i\}$ must be infinite;
for example, let $|\{u_i\}|=\infty$.

It follows from $\gamma(u_i,v_i)\ne 0$ that for each $i$ there is $z_i\in [u_i,v_i]$ such
that $\alpha(u_i,z_i)\ne 0\ne \beta(z_i,v_i)$. Since $\alpha\in FI(P)$ and
$[u_i,z_i]\subset [x,y]$, we have $u_i=z_i$ for an infinite number of indexes. But then
$\beta(u_i,v_i)\ne 0$ for this set of indexes, what is impossible, since $u_i\ne v_i$.
$\blacksquare$

\medskip

{\bf Example.} Let $\mathbb{N}$ be the set of positive integers with
the natural order, $\overline{\mathbb{N}}$ its isomorphic copy and
$P=\mathbb{N}\oplus\overline{\mathbb{N}}$ their order sum (i.\,e.
$a<\overline{b}$ for each $a\in \mathbb{N}$, $\overline{b}\in
\overline{\mathbb{N}}$).

The algebra $FI(\overline{\mathbb{N}})$ embeds into $FI(P)$; it is sufficient for this to
extend each series from $FI(\overline{\mathbb{N}})$ by zero values. This is not the case
for $FI(\mathbb{N})$: since each segment of the form $[a,\overline{b}]$ ($a\in
\mathbb{N}$, $\overline{b}\in \overline{\mathbb{N}}$) is infinite, each series from
$FI(P)$ is nonzero only at a finite number of the segments from $\mathbb{N}$. Therefore
$FI(P)$ contains not $FI(\mathbb{N})$, but the algebra of {\it finite} formal sums of the
segments from $\mathbb{N}$.

\section{Properties of finitary algebras}\label{sec-2}

Denote by $\delta$ the identity of the algebra $FI(P)$. Then $\delta(x,y)=\delta_{xy}$
where $\delta_{xy}$ is the Kronecker delta.

\begin{thrm}\label{thrm-2}
A series $\alpha\in FI(P)$ is invertible iff $\alpha(x,x)\ne 0$
for each $x\in P$. Moreover $\alpha^{-1}\in FI(P)$.
\end{thrm}

{\bf Proof.} \emph{Necessity.} Let $\alpha\beta=\delta$. Then
$\alpha(x,x)\beta(x,x)=\delta(x,x)=1$, so $\alpha(x,x)\ne 0$.

\emph{Sufficiency.} Let $\alpha(x,x)\ne 0$ for all $x\in P$ and $[u,v]$ be a segment from
$P$, $u\ne v$. A series $\beta$, which is inverse for $\alpha$, exists if
$\beta(x,x)=\alpha(x,x)^{-1}$ and
\begin{equation}\label{eq-2}
\beta(u,v)=-\alpha(u,u)^{-1}\sum_{u< x\le v}\alpha(u,x)\beta(x,v)
\end{equation}
(note that the sum on the right-hand side is defined since $\alpha(u,x)$ is different
from zero only for a finite number of elements $x\in [u,v]$). We prove that a solution of
the equation~\eqref{eq-2} exists and can be computed recursively in a finite number of
steps.

Denote by $C_{\alpha}(u,v)$ the number of subsegments
$[x,y]\subseteq [u,v]$ such that $x\ne y$ and $\alpha(x,y)\ne 0$. By
the definition of finitary series $C_{\alpha}(u,v)$ is finite. We
shall prove our assertion by induction on $C_{\alpha}(u,v)$.

If $C_{\alpha}(u,v)=0$ then $\beta(u,v)= 0$. If $C_{\alpha}(u,v)=1$ and $u<x_0\le v$,
$\alpha(u,x_0)\ne 0$ then
\begin{eqnarray}
\beta(u,v)&=&-\alpha(u,u)^{-1}\alpha(u,x_0)\beta(x_0,v)\nonumber\\
&=& \left\{
\begin{array}{l}
-\alpha(u,u)^{-1}\alpha(u,x_0)\alpha(v,v)^{-1}, \ \ \  \mbox{\rm if}\ x_0=v,\\
0, \ \ \ \ \ \ \ \ \ \ \ \mbox{\rm if}\ x_0\ne v \ (\mbox{\rm
since} \ C_{\alpha}(x_0,v)=0).
\end {array} \right.\nonumber
\end {eqnarray}

Now suppose that $\beta(x,y)$ is defined for all $x,y$ such that $C_{\alpha}(x,y)<n$. Let
$C_{\alpha}(u,v)=n$. If $\alpha(u,x)$ is nonzero for some $x$, $u<x\le v$, then
$C_{\alpha}(x,v)\le n-1$ and, by the induction hypothesis, $\beta(x,v)$ is defined. Thus,
every summand from the right-hand side of~\eqref{eq-2} is defined, and since the sum is
finite, $\beta(u,v)$ is also defined.

Suppose that $\alpha^{-1}\not\in FI(P)$. Then we can find a segment $[x,y]$ and an
infinite number of subsegments $[u_i,v_i]\subset [x,y]$, $u_i\ne v_i$, for which
$\alpha^{-1}(u_i,v_i)\ne 0$. From the equalities
$\alpha^{-1}\alpha=\alpha\alpha^{-1}=\delta$ we have:
\begin{equation}\label{eq-3}
\alpha^{-1}(u_i,v_i)=-\alpha(v_i,v_i)^{-1}\sum_{u_i\le z_i<v_i}
\alpha^{-1}(u_i,z_i)\alpha(z_i,v_i),
\end{equation}
\begin{equation}\label{eq-4}
\alpha^{-1}(u_i,v_i)=-\alpha(u_i,u_i)^{-1}\sum_{u_i< z_i\le v_i}
\alpha(u_i,z_i)\alpha^{-1}(z_i,v_i).
\end{equation}

It follows from~\eqref{eq-3} that for each $i$ there is $z_i$ such that
$\alpha(z_i,v_i)\ne 0$. But the number of such segments is finite, therefore
$|\{v_i\}|<\infty$. Similarly $|\{u_i\}|<\infty$ from~\eqref{eq-4}, a contradiction.
Hence $\alpha^{-1}\in FI(P)$. $\blacksquare$

\begin{cor}\label{cor-1} Left, right and two-sided invertibilities in $FI(P)$ coincide.
\end{cor}

{\bf Proof.} If, for example, $\alpha$ is right invertible, then, as was shown in the
proof of necessity, $\alpha(x,x)\ne 0$ for all $x\in P$. Now it follows from the theorem,
that $\alpha$ is two-sided invertible. $\blacksquare$

\bigskip

In what follows ${\rm Rad}\,A$ denotes the Jacobson radical of
algebra $A$.

\begin{cor}\label{cor-2} $\alpha\in {\rm Rad}\,FI(P)$ iff
$\alpha(x,x)=0$ for all $x\in P$.
\end{cor}

{\bf Proof.} As in~\cite{d-s-r}, the proof follows
from~\cite{jac}, Proposition 1.6.1. $\blacksquare$

\begin{cor}\label{cor-3} The factor algebra $FI(P)/{\rm Rad}\,FI(P)$ is commutative.
$\blacksquare$
\end{cor}

The following proposition describes the idempotents of a finitary
algebra.

We call an element $\alpha\in FI(P)$ {\it diagonal} if
$\alpha(x,y)=0$ for $x \ne y$. It is obvious that a diagonal element
$\alpha$ is idempotent iff $\alpha(x,x)$ is equal to 0 or 1 for all
$x\in P$.

\begin{thrm}\label{thrm-3}
Each idempotent $\alpha \in FI(P)$ is conjugate to the diagonal
idempotent $\varepsilon$, such that $\varepsilon(x,x)=\alpha(x,x)$
for all $x\in P$.
\end{thrm}
{\bf Proof.} According to the corollary~\ref{cor-2}, $\rho=\alpha-\varepsilon \in{\rm
Rad}\,FI(P)$. Since $\alpha^2=\alpha$, we have:
\begin{equation}\label{eq-5}
\varepsilon\rho+\rho\varepsilon=\rho-\rho^2
\end{equation}
Multiplying this equality by $\varepsilon$ on the left, we obtain:
\begin{equation}\label{eq-6}
\varepsilon\rho\varepsilon+\varepsilon\rho^2=0
\end{equation}
Set $\beta=\delta+(2\varepsilon-\delta)\rho$ where $\delta$ is the identity of $FI(P)$.
Since $\rho\in{\rm Rad}FI(P)$, the series $\beta$ is invertible by theorem~\ref{thrm-2}.

By~\eqref{eq-6}, we have:
\begin{eqnarray*}
\beta\alpha&=&(\delta+2\varepsilon\rho-\rho)(\varepsilon+\rho)=
\varepsilon-\rho\varepsilon+\rho-\rho^2,\\
\varepsilon\beta&=&\varepsilon(\delta+2\varepsilon\rho-\rho)=\varepsilon+\varepsilon\rho.
\end{eqnarray*}
From~\eqref{eq-5} we obtain: $\beta\alpha=\varepsilon\beta$. Hence
$\alpha=\beta^{-1}\varepsilon\beta$. $\blacksquare$

\bigskip

An element $\alpha$ of the algebra $A$ is called {\it regular} if there is $\chi \in A$,
such that $\alpha\chi\alpha=\alpha$.

\begin{thrm}\label{thrm-4}
For each regular $\alpha\in FI(P)$ there are a diagonal idempotent
$\varepsilon \in FI(P)$ and invertible elements $\beta, \gamma \in
FI(P)$, such that $\alpha=\beta\varepsilon\gamma$.
\end{thrm}
{\bf Proof.} The sufficiency is obvious --- we can set
$\chi=\gamma^{-1}\varepsilon\beta^{-1}$. Let us prove the
necessity.

Evidently, $\alpha\chi$ and $\chi\alpha$ are idempotents. By the
theorem~\ref{thrm-3}, there is an invertible $\eta \in FI(P)$, such
that $\alpha\chi=\eta^{-1}\varepsilon\eta$, where $\varepsilon$ is a
diagonal idempotent, and, moreover,
$\varepsilon(x,x)=\alpha(x,x)\chi(x,x)$. In consequence of the
regularity, the last equality is equivalent to the statement
\begin{equation}\label{eq-7}
\varepsilon(x,x)=0 \Longleftrightarrow \alpha(x,x)=0.
\end{equation}

Similarly, $\chi\alpha=\gamma^{-1}\varepsilon_1\gamma$ for some invertible $\gamma \in
FI(P)$ and diagonal idempotent $\varepsilon_1 \in FI(P)$. In fact $\varepsilon_1
=\varepsilon$ since~\eqref{eq-7} holds for $\varepsilon_1$ too.

From the regularity of $\alpha$ we have:
\[
\alpha=\alpha\chi\alpha=\alpha\chi\alpha\chi\alpha=
\eta^{-1}\varepsilon\eta\alpha\gamma^{-1}\varepsilon\gamma.
\]

Since $\eta$ and $\gamma$ are invertible, it follows from~\eqref{eq-7} that
$$
\eta\alpha\gamma^{-1}(x,x)=0 \Longleftrightarrow
\varepsilon(x,x)=0.
$$
Therefore the element $\eta\alpha\gamma^{-1}$ can be rewritten as
$\eta\alpha\gamma^{-1}=\eta_1\varepsilon+\rho$ where $\eta_1$ is
diagonal invertible and $\rho\in {\rm Rad}\,FI(P)$. Since diagonal
elements commute, $\varepsilon\eta\alpha\gamma^{-1}\varepsilon=
\eta_1\varepsilon+\varepsilon\rho\varepsilon=(\eta_1+\varepsilon\rho)\varepsilon$,
and, moreover, $\eta_1+\varepsilon\rho$ is invertible by the
theorem~\ref{thrm-2}. Thus,
$$
\alpha=\eta^{-1}\varepsilon\eta\alpha\gamma^{-1}\varepsilon\gamma=
\eta^{-1}(\eta_1+\varepsilon\rho)\varepsilon\gamma.\ \blacksquare
$$

In conclusion we consider a property of elements, which is intermediate between the
invertibility and the regularity.

Let $\alpha$ be a regular element, $\alpha\chi\alpha=\alpha$. Then the element
$\alpha^\ast=\chi\alpha\chi$ satisfies the equations:
\begin{equation}\label{eq-8}
\alpha\alpha^\ast\alpha=\alpha, \ \
\alpha^\ast\alpha\alpha^\ast=\alpha^\ast.
\end{equation}

We call $\alpha$ {\it superregular} if there exists only one $\alpha^\ast$ for which the
equations~\eqref{eq-8} are fulfilled. For instance, invertible and zero elements are
superregular. It turns out that in incidence algebras superregular elements can be
described with the help of invertible ones:

\begin{cor}\label{cor-4}
Let $P=\bigcup\limits_{i\in I}P_i$ be decomposition of $P$ into a disjoint union of
connected components and $U(P_i)$ be the group of invertible elements of the algebra
$FI(P_i)$. Then the set of superregular elements coincides with the sum of semigroups
$\bigoplus\limits_{i\in I}U(P_i)^0$.
\end{cor}
{\bf Proof.} It is easy to see that $FI(P)=\bigoplus\limits_{i\in
I}FI(P_i)$. Obviously, each superregular element $\alpha\in FI(P)$
can be represented as the sequence $(\alpha_i)_{i\in I}$ of the
superregular elements $\alpha_i\in FI(P_i)$ and, moreover,
$\alpha^\ast=(\alpha_i^\ast)_{i\in I}$. The statement will be proven
if we shall show that $\alpha_i\in U(P_i)^0$.

Suppose the contrary. First note that the multiplication on the left or on the right by
an invertible element preserves superregularity. Therefore, by the theorem~\ref{thrm-4},
we can consider $\alpha_i$ to be a diagonal idempotent. By assumption, there is a pair
$x,y\in P_i$, such that $\alpha_i(x,x)\ne \alpha_i(y,y)$, and, moreover, we can choose it
in such a way, that $x<y$ since $P_i$ is connected. Define $\alpha^\ast_i=\alpha_i+\rho$,
where
\begin{equation*}
\rho(u,v)=\begin{cases}
   1, & \mbox{if $u=x,v=y$},\\
   0, & \mbox{otherwise}.
   \end{cases}
\end{equation*}
By the direct checking, we make sure, that the equalities~\eqref{eq-8} hold for
$\alpha_i$, and, since $\alpha_i$ is a superregular idempotent, $\alpha^\ast_i=\alpha_i$.
This contradicts the definition of $\rho$. $\blacksquare$

\section{The isomorphism problem}\label{sec-3}

It is well known~\cite{d-s-r} that the isomorphism problem for
locally finite posets is solved positively: if $P$ and $Q$ are
locally finite and $I(P) \cong I(Q)$, then $P \cong Q$ (in the case
$K$ is a ring, an answer can be negative, see, for
example,~\cite{vos}). In this section we consider the isomorphism
problem for finitary series.

Recall, that an idempotent $\alpha\ne 0$ is {\it primitive}, if
$\alpha\varepsilon=\varepsilon\alpha=\alpha$ for some idempotent $\varepsilon$ implies
that $\varepsilon$ is equal to 0 or $\alpha$.

We will need the idempotents of special form $\delta_x^P$, defined
for an arbitrary element $x\in P$ as follows:
\begin{eqnarray}
\delta_x^P(u,v) &=& \begin{cases}
1, & \mbox{if $u = v = x$},\\
0, & \mbox{if $u \neq x \text{ or } v \neq x$}.
\end{cases}\label{eq-9}
\end{eqnarray}

\begin{lem}\label{lem-1}
An idempotent $\alpha$ is primitive iff it is conjugate to
$\delta_x^P$ for some $x\in P$.
\end{lem}

{\bf Proof.} The conjugation, being an automorphism, preserves the primitivity. So by the
theorem~\ref{thrm-3}, it is sufficient to prove that the diagonal primitive idempotents
are nothing but $\delta_x^P$.

\medskip

1) We first show that $\delta_x^P$ is primitive. Indeed, let $\alpha\delta_x^P =\alpha$
holds for some idempotent $\alpha$. Then for each segment $[u,v] \subset P$ we obtain:
$$
\alpha(u,v)\delta_x^P(v,v)-\alpha(u,v)=0,
$$
i.\,e. $\alpha(u,v)=0$ if $v\neq x$.

Similarly $\delta_x^P\alpha = \alpha$ is equivalent to $\alpha(u,v)=0$ for $u\neq x$. If
$\alpha(x,x)=0$ then we get $\alpha \equiv 0$. On the other hand, if $\alpha(x,x)=1$ then
$\alpha = \delta_x^P$.

\medskip

2) Let $\varepsilon$ be a diagonal primitive idempotent.

Since $\varepsilon \neq 0$, there is such $x \in P$ that $\varepsilon(x,x) = 1$. Consider
the idempotent $\delta_x^P$. Obviously, the equalities $\varepsilon\delta_x^P =
\delta_x^P = \delta_x^P\varepsilon$ hold for it. Therefore $\delta_x^P = \varepsilon$,
since $\varepsilon$ is primitive and $\delta_x^P \neq 0$. $\blacksquare$

\bigskip

It is easy to see that if $\delta_x^P$ and $\delta_y^P$ are conjugate then $x=y$.

\begin{thrm}\label{thrm-5}
Let $P$ and $Q$ be arbitrary posets. Then
\[
FI(P) \cong FI(Q) \Longrightarrow P \cong Q
\]
\end{thrm}

{\bf Proof.} Let $\Phi:FI(P)\rightarrow FI(Q)$ be an isomorphism. For each $x\in P$ the
image $\Phi(\delta_x^P)$ of $\delta_x^P$ is primitive and, by lemma~\ref{lem-1}, is
conjugate to the idempotent $\delta_y^Q$ for some $y\in Q$. According to the remark
before theorem, the element $y$ is defined uniquely. Thus, $\Phi$ generates the bijection
$\varphi: P \rightarrow Q$, such that $\Phi(\delta_x^P)$ is conjugate to
$\delta_{\varphi(x)}^Q$.

Let us prove that $\varphi$ preserves the order.

It is easy to see that for all $x,y \in P$
\[
x \leq y \Longleftrightarrow \delta_x^P FI(P)\delta_y^P \neq 0.
\]

Similarly for all $u,v \in Q$
\[
u \leq v \Longleftrightarrow \delta_u^QFI(Q)\delta_v^Q \neq 0.
\]
Let $x \leq y$ and $\Phi(\delta_x^P)=\beta_1\delta_{\varphi(x)}^Q\beta_1^{-1}$,
$\Phi(\delta_y^P)=\beta_2\delta_{\varphi(y)}^Q\beta_2^{-1}$ for some invertible elements
$\beta_1$ and $\beta_2$. From the bijectivity of $\Phi$ we obtain
\[
\Phi(\delta_x^P)FI(Q)\Phi(\delta_y^P) = \Phi(\delta_x^P
FI(P)\delta_y^P) \neq 0.
\]
Hence
\[
\delta_{\varphi(x)}^QFI(Q)\delta_{\varphi(y)}^Q =
\beta_1^{-1}\Phi(\delta_x^P)\beta_1FI(Q)\beta_2^{-1}\Phi(\delta_y^P) \beta_2 =
\beta_1^{-1}\Phi(\delta_x^P)FI(Q)\Phi(\delta_y^P)\beta_2
\]
since $\beta_1$ and $\beta_2$ are invertible. Therefore
$\varphi(x)\leq \varphi(y)$. $\blacksquare$

\begin{cor}\label{cor-5}
Let $P$ be not locally finite. Then there is no locally finite poset $Q$, such that
$FI(P)\cong I(Q)$. $\blacksquare$
\end{cor}

\end{document}